\documentclass{amsproc}
\usepackage{eurosym}
\usepackage{amssymb}
\usepackage{amsfonts}
\usepackage{cmtt}

\theoremstyle{plain}
\newtheorem{theorem}{Theorem}

\theoremstyle{remark}

\begin{document}

\title{Matroid multiple cyclic exchange property}

\author[Micha{\l} Laso\'{n}]{Micha{\l} Laso\'{n}}

\dedicatory{\upshape
Institute of Mathematics of the Polish Academy of Sciences,\\ ul.\'{S}niadeckich 8, 00-656 Warszawa, Poland\\ \textmtt{michalason@gmail.com}}

\thanks{Research supported by Polish National Science Centre grant no. 2015/19/D/ST1/01180. The paper was completed during author's stay at Freie Universit\"at Berlin in the frame of Polish Ministry ``Mobilno\'s\'c Plus'' program.}

\begin{abstract}
We prove a new exchange property for bases of a matroid that generalizes the multiple symmetric exchange property. For every bases $B_1,\dots,B_k$ of a matroid and a subset $A_1\subset B_1$ there exist subsets $A_2\subset B_2,\dots,A_k\subset B_k$ such that all sets $(B_i\setminus A_i)\cup A_{i-1}$ achieved by a cyclic shift of $A_i$'s by one are bases. 
\end{abstract}

\maketitle


Base exchange properties is a central theme in matroid theory. A non-empty family $\mathfrak{B}$ of subsets of a finite set $E$ forms a set of bases of a matroid, just from the definition, if it satisfies the \emph{exchange property}. That is, if for every $B_1,B_2\in\mathfrak{B}$ and $e_1\in B_1\setminus B_2$ there exists $e_2\in B_2\setminus B_1$, such that $(B_1\setminus e_1)\cup e_2\in\mathfrak{B}$ (we refer the reader to \cite{Ox92} for a background of matroid theory). 

In this case a stronger property holds. For every bases $B_1,B_2$ of a matroid and an element $e_1\in B_1\setminus B_2$ there exists an element $e_2\in B_2\setminus B_1$, such that both sets $(B_1\setminus e_1)\cup e_2$ and $(B_2\setminus e_2)\cup e_1$ are bases. It is the \emph{symmetric exchange property} discovered by Brualdi \cite{Br69} (White's conjecture on the toric ideal of a matroid is based on this property, see \cite{LaMi14,La16,Wh80}).

Even more is true matroid bases -- one can exchange symmetrically not only single elements, but also subsets. The \emph{multiple symmetric exchange property} asserts that for every bases $B_1,B_2$ of a matroid and a subset $A_1\subset B_1$ there exists a subset $A_2\subset B_2$, such that both sets $(B_1\setminus A_1)\cup A_2$ and $(B_2\setminus A_2)\cup A_1$ are bases. It was proved by Green \cite{Gr73}, and soon after by Woodall \cite{Wo74} (see \cite{Ku86} for more exchange properties, and \cite{La15} for easy proofs of them).

In this note we formulate and prove a new exchange property that generalizes the multiple symmetric exchange property.

\begin{theorem}[Multiple cyclic exchange property]\label{TheoremMain}
For every bases $B_1,\dots,B_k$ of a matroid and a subset $A_1\subset B_1$ there exist subsets $A_2\subset B_2,\dots,A_k\subset B_k$ such that all sets $(B_i\setminus A_i)\cup A_{i-1}$ achieved by a cyclic shift of $A_i$'s by one are bases. 
\end{theorem}

Notice that for $k=2$ we get exactly the multiple symmetric exchange property.

\begin{proof}
The proof uses a link between list coloring of a matroids and base exchange properties, see \cite{La15}. We consider a certain assignment of lists of colors to elements of the ground set, and show that there exists a proper (elements of the same color form an independent set) coloring from these lists. This implies a desired exchange property.

Let $M$ be a matroid on a ground set $E$ with rank function $r$. Suppose $B_1,\dots,B_k$ are bases of $M$, and $A_1\subset B_1$. We can assume that bases $B_i$ are pairwise disjoint (since if they are not, then we can introduce parallel elements).

Let $L$ assign list $\{1\}$ to elements of $B_1\setminus A_1$, list $\{2\}$ to elements of $A_1$, list $\{i,i+1\}$ to elements of $B_i$ for $i=2,\dots,k-1$, and list $\{1,k\}$ to elements of $B_k$. For every color $i$ consider the set of elements $C_i$ that have $i$ on their list
$$C_1=B_k\cup (B_1\setminus A_1),C_2=A_1\cup B_2,C_3=B_2\cup B_3,\dots,C_k=B_{k-1}\cup B_k,$$   
and the matroid $M_i$ with rank function $r_i$, achieved by restricting $M$ to the set $C_i$. 

Observe that for every $A\subset B_1\cup\cdots\cup B_k$ inequality $r_1(A)+\dots+r_k(A)\geq\lvert A\rvert$ holds. Indeed, 
$r_1(A)+\dots+r_k(A)=$
$$r(A\cap (B_k\cup (B_1\setminus A_1)))+r(A\cap (A_1\cup B_2))+r(A\cap (B_2\cup B_3))+\dots+r(A\cap (B_{k-1}\cup B_k))$$
using submodularity of rank function to the first and the last summand we obtain
$$\geq r(A\cap (B_{k-1}\cup B_k\cup (B_1\setminus A_1)))+r(A\cap (A_1\cup B_2))+r(A\cap (B_2\cup B_3))+\dots+r(A\cap B_k)\geq$$
$$\geq r(A\cap (B_{k-1}\cup (B_1\setminus A_1)))+r(A\cap (A_1\cup B_2))+r(A\cap (B_2\cup B_3))+\dots+r(A\cap B_k)$$
and inductively
$$\geq r(A\cap (B_2\cup (B_1\setminus A_1)))+r(A\cap (A_1\cup B_2))+r(A\cap B_3)+\dots+r(A\cap B_k)\geq$$
$$\geq r(A\cap B_1)+r(A\cap B_2)+r(A\cap B_3)+\dots+r(A\cap B_k)=\lvert A\cap B_1\rvert+\dots+\lvert A\cap B_k\lvert=\lvert A\rvert.$$

By the matroid union theorem there exist sets $D_1,\dots,D_k$, each $D_i$ independent in $M_i$, such that $D_1\cup\cdots\cup D_k=B_1\cup\cdots\cup B_k$. Observe that $A_1=D_2\cap B_1$.
Define $A_i=D_{i+1}\cap B_i$ for $i=2,\dots,k$. Then $(B_i\setminus A_i)\cup A_{i-1}=D_i$, so it is a basis.
\end{proof}

For graphic matroids the single element cyclic exchange property (when $\lvert A_1\rvert$=1) was firstly proved by Bart{\l}omiej Bosek \cite{Bo16} using cycle axioms.

Theorem \ref{TheoremMain} is optimal in a sense that for every $k\geq 3$ there exists a matroid of rank $3$ and $k$ bases $B_1,\dots,B_k$ for which one can not have additionally that all sets $(B_i\setminus A_i)\cup A_{i-2}$ achieved by a cyclic shift of $A_i$'s by two are bases.


\end{document}